\documentclass{amsart}

\usepackage{a4wide}
\usepackage{amsmath}
\usepackage{enumerate}

\usepackage{amsmath,amsthm}
\usepackage{amsfonts}
\usepackage{amssymb}
\usepackage{graphicx, color, soul}
\usepackage[T1]{fontenc}
\usepackage{verbatim}

\newtheorem*{theorem*}{Theorem}
\newtheorem*{cor*}{Corollary}

\theoremstyle{plain}

\newtheorem{theorem}{Theorem}[section]

\newtheorem{prop}[theorem]{Proposition}

\newtheorem{lem}[theorem]{Lemma}

{\rm}{\rm}

\theoremstyle{definition}

\begin{document}
	
\title{On holonomy groups of K-contact sub-pseudo-Riemannian manifolds}
	\author{E.A. Kokin}

	\begin{abstract} 
		This article investigates the holonomy groups of K-contact sub-pseudo-Riemannian manifolds. The primary result is a proof that the horizontal holonomy group either coincides with the adapted holonomy group or acts as its normal subgroup of codimension one. The theory is adapted for metrics of indefinite signature, bypassing the problem of subspace degeneracy that previously prevented the use of established orthogonal decomposition methods. It is established that, in the sub-Lorentzian case, the adapted holonomy group corresponds to the holonomy group of a certain Lorentzian manifold. This work also provides a complete classification of codimension-one ideals for Lorentzian holonomy algebras and presents specific examples of structures based on Cahen-Wallach spaces and Kähler manifolds.
		
		
	\end{abstract}

	\maketitle

	\section{Introduction}
	
	The holonomy group is an important invariant of a connection on a smooth manifold. An essential role is played by the classification of holonomy groups of Riemannian manifolds \cite{Besse,bryant99,Joyce07}. There exists a classification of holonomy groups of Lorentzian manifolds \cite{BBI,Gal06,
		Leistner}, as well as particular results on the holonomy groups of pseudo-Riemannian manifolds, see, for example, \cite{FK,GalLei10,
		GalpsKahl}.
	
	Sub-Riemannian manifolds are important from the perspective of geometric control theory and rolling problems \cite{Agr19,Montg}. In recent works, interest has also been shown in sub-Lorentzian structures \cite{CMV08,Markina16,Molina2025}. 	
	Motivated by recent papers \cite{GalHolK,GLL25} on the holonomy of contact sub-Riemannian manifolds, in the present article we begin the study of holonomy groups of contact sub-Lorentzian manifolds and, more generally, sub-pseudo-Riemannian manifolds.
	
	A sub-pseudo-Riemannian manifold $(M,\theta,g)$   is a smooth manifold $M$ with a fixed contact form $\theta$ and a field $g$ of non-degenerate symmetric bilinear forms on the contact distribution $D=\ker\theta$. 	
	Such a structure defines a horizontal connection $\nabla^g$, which specifies the parallel transport of vectors tangent to $D$ along curves tangent to $D$, and we obtain the corresponding holonomy group $\operatorname{Hol}_x(\nabla^g)$, $x\in M$. 	
	There also exists an adapted connection $\nabla^{\boldsymbol{\tau}}$, which extends the connection $\nabla^g$ to a connection on the vector bundle $D$ over $M$. Let $\operatorname{Hol}_x(\nabla^{\boldsymbol{\tau}})$ denote the corresponding holonomy group at point $x\in M$. 	
	A sub-pseudo-Riemannian manifold $(M,\theta,g)$ is called $K$-contact if $L_\xi g=0$, where $L_\xi$   is the Lie derivative along the Reeb field $\xi$.	
	First of all, we prove the following theorem.
	
	\begin{theorem}\label{ThpsRiem}
		Let $(M,\theta,g)$ be a $K$-contact sub-pseudo-Riemannian manifold. Let $x\in M$. Then one of the following conditions holds:
		\begin{itemize}
			\item[1.] $\operatorname{Hol}_x(\nabla^g)=\operatorname{Hol}_x(\nabla^{\boldsymbol{\tau}})$;
			\item[2.] $\operatorname{Hol}_x(\nabla^g)\subset\operatorname{Hol}_x(\nabla^{\boldsymbol{\tau}})$ is a normal subgroup of codimension one.
		\end{itemize}
	\end{theorem}
	
	Theorem \ref{ThpsRiem} was proven in \cite{GalHolK} for $K$-contact sub-Riemannian manifolds under the assumption that the Reeb field $\xi$ is complete. Later, this theorem was proven in \cite{GLL25} for $K$-contact sub-Riemannian manifolds with an arbitrary $\xi$. 
	The proof from \cite{GLL25} cannot be extended to sub-pseudo-Riemannian manifolds because the authors of \cite{GLL25} considered the orthogonal complement $\mathfrak{hol}_x(\nabla^g)^\bot$ of the holonomy algebra $\mathfrak{hol}_x(\nabla^g)$ in the holonomy algebra $\mathfrak{hol}_x(\nabla^{\boldsymbol{\tau}})$ and the orthogonal decomposition
	$$\mathfrak{hol}_x(\nabla^{\boldsymbol{\tau}})=\mathfrak{hol}_x(\nabla^g)\oplus \mathfrak{hol}_x(\nabla^g)^\bot.$$
	In the case of sub-pseudo-Riemannian manifolds, it is possible that $\mathfrak{hol}_x(\nabla^g)\subset\mathfrak{hol}_x(\nabla^{\boldsymbol{\tau}})$ is a degenerate subspace, and then the aforementioned orthogonal decomposition cannot be considered. Our proof of Theorem \ref{ThpsRiem} uses some ideas from \cite{GalHolK}, as well as new ideas.
	
	Next, we consider the case of $K$-contact sub-Lorentzian spaces. First of all, \cite[Prop. 3]{GalHolK} can be directly extended to
	
	\begin{prop}\label{PropAdaHol}
		Let $(M,\theta,g)$ be a $K$-contact sub-Lorentzian manifold of dimension $n=2m+1$. Then the identity component of the adapted holonomy group $\operatorname{Hol}^0(\nabla^{\boldsymbol{\tau}})\subset \operatorname{SO}(1,2m-1)$
		is the holonomy group of a Lorentzian manifold.
	\end{prop}
	
	We see that horizontal holonomy algebras are exhausted by Lorentzian holonomy algebras and codimension-one ideals of Lorentzian holonomy algebras. We list these algebras in Section \ref{SecLor}. We also consider examples.
	
	\section{Preliminaries}
	
	We use the exposition from \cite{GalHolK,GLL25}, which relies on results from \cite{CGJK,FGR} and which generalizes directly from the case of contact sub-Riemannian manifolds to the case of contact sub-pseudo-Riemannian manifolds.
	
	Let $(M,\theta,g)$   be a {\bf contact sub-pseudo-Riemannian manifold}. Here $M$   is a smooth manifold of dimension $n=2m+1\ge 5$ with a fixed contact form $\theta$ and a field $g$ of non-degenerate symmetric bilinear forms on the contact distribution $D=\ker\theta$. 
	
	Let $\xi$ denote the corresponding {\bf Reeb field}, uniquely defined by two conditions:
	$$\theta(\xi)=1,\quad \iota_\xi (d\theta)=0.$$
	This gives the direct sum
	$$TM=\left<\xi\right>\oplus D$$
	and projections $\pi:TM\to D$, $\pi':TM\to \left<\xi\right>$.
	
	The metric $g$ defines the {\bf horizontal connection} (an analogue of the Levi-Civita connection)
	$$\nabla^g:\Gamma(D)\times \Gamma(D)\to\Gamma(D)$$
	by the equality
	\begin{equation}\label{SchoutenCon}
	2g(\nabla^g_X Y,Z)= Xg(Y,Z) +Yg(Z,X) - Zg(X,Y)+ g(\pi[X,Y],Z)- g(\pi[Y,Z],X)- g(\pi[X,Z],Y),
	\end{equation} 
	for all $X,Y,Z\in\Gamma(D)$.
	
	The Schouten curvature tensor of the connection $\nabla^g$ is defined as
	\begin{equation}\label{Schouten}
	R^g(X,Y)Z=\nabla^g_X\nabla^g_Y Z-\nabla^g_Y\nabla^g_X Z-\nabla^g_{\pi[X,Y]}Z-\pi\left([\pi' ([X,Y]),Z]\right),\quad \text{for all } X,Y,Z\in D.
	\end{equation}
	
	The horizontal connection $\nabla^g$ defines the parallel transport $\tau^g_\gamma$ of vectors tangent to $D$ along curves $\gamma$ tangent to $D$, and the corresponding holonomy group $\operatorname{Hol}_x(\nabla^g)$, $x\in M$, is called the {\bf horizontal holonomy group}.
	
	Any endomorphism $N:D\to D$ defines an {\bf extended connection} \cite{Gal3,Gal1} $$\nabla^N:\Gamma(TM)\times \Gamma(D)\to\Gamma(D)$$ by the equalities
	$$\nabla^N_XY=\nabla^g_XY,\quad \nabla^N_\xi Y=[\xi,Y]+NY,\quad X,Y\in\Gamma(D).$$
	The connection $\nabla^N$ is a connection in the vector bundle $D$ over the manifold $M$. Let $R^N$, $\tau^N$ and $\operatorname{Hol}_x(\nabla^N)$ denote the corresponding curvature tensor, parallel transport and holonomy group.
	
	We will use the symmetric endomorphism ${\boldsymbol{\tau}}$, defined as
	$$g({\boldsymbol{\tau}} X,Y)=\frac{1}{2}L_\xi g (X,Y),\quad X,Y\in\Gamma(D).$$
	We consider two extended connections. The first is the {\bf adapted connection} $\nabla^{\boldsymbol{\tau}}$, defined by the equalities 
	$$\nabla^{\boldsymbol{\tau}}_XY=\nabla^g_XY,\quad \nabla^{\boldsymbol{\tau}}_\xi Y=[\xi,Y]+{\boldsymbol{\tau}} Y,\quad X,Y\in\Gamma(D).$$
	
	The second is the {\bf Wagner connection} $\nabla^W$ \cite{Gal3,Wagner}. The corresponding endomorphism 
	$$N^W=\frac{1}{4m}R((d\theta)^{-1}),$$ 
	where $(d\theta)^{-1}$   is the tensor field of type $(2,0)$, inverse to the tensor field $d\theta$.
	
	The following theorem is important, following from the results of \cite{CGJK,FGR}, see \cite{GLL25}.
	
	\begin{theorem}\label{thholW} 
		Let $(M,\theta,g)$ be a contact sub-pseudo-Riemannian manifold and $x\in M$. Then the horizontal holonomy group $\operatorname{Hol}_x(\nabla^g)$ coincides with the holonomy group $\operatorname{Hol}_x(\nabla^W)$ of the Wagner connection.
	\end{theorem}
	
	A sub-pseudo-Riemannian manifold $(M,\theta,g)$ is called {\bf $K$-contact}, if
	$${\boldsymbol{\tau}}=0.$$ 
	This means that $L_\xi g=0$. 
	In this case $N^W$ is skew-symmetric and we denote it by $C$. 
	In this case it holds 
	$$\nabla^W_\xi=\nabla^{\boldsymbol{\tau}}_\xi+C.$$ 
	Finally, the following equalities hold 
	$$R^{\boldsymbol{\tau}}(X,Y)=R^g(X,Y),\quad R^{\boldsymbol{\tau}}(\xi,X)=0,\quad X,Y\in\Gamma(D).$$
	
	\section{Proof of theorem \ref{ThpsRiem}}
	
	Since ${\boldsymbol{\tau}}=0$, then $N^W=C$ is a skew-symmetric operator.
	
	\begin{lem}\label{lemxi1} Let $\lambda:[0,r]\to M$ be an integral curve of the vector field $\xi$. Let $x=\lambda(0)$ and $y=\lambda(r)$. Then the following is true
		$$\tau_\lambda^{\boldsymbol{\tau}}=\tau_\lambda^W\circ e^{rC_x}=e^{rC_y}\circ \tau_\lambda^W.$$
	\end{lem}
	
	{\bf Proof.} Let $X(s)$ be a $\nabla^{\boldsymbol{\tau}}$-parallel section of $D$ along the curve $\lambda$.
	Then, since $\dot\lambda(s)=\xi_{\lambda(s)}$ and $\nabla^{\boldsymbol{\tau}}_\xi C=L_\xi C=0$, we have
	\begin{multline*}
	\nabla_{\dot\lambda(s)}^W\left(e^{-sC_{\lambda(s)}}X(s) \right)=\left(\nabla_{\dot\lambda(s)}^{\boldsymbol{\tau}}+C_{\lambda(s)}\right)\left(e^{-sC_{\lambda(s)}}X(s) \right)\\  -C_{\lambda(s)}e^{-sC_{\lambda(s)}}X(s)   + C_{\lambda(s)}e^{-sC_{\lambda(s)}}X(s) =0.\end{multline*}
	This shows that $$\tau^W_\lambda X(0)=e^{-rC_y}X(r)=
	e^{-rC_y}\tau^{\boldsymbol{\tau}}_\lambda X(0),$$
	whence it follows that $$\tau^{\boldsymbol{\tau}}_\lambda=e^{rC_y}\circ 
	\tau^W_\lambda.$$ Consider the curve $\lambda_1(s)=\lambda(r-s)$, $s\in[0,r]$. Let $X(s)$ be a $\nabla^{\boldsymbol{\tau}}$-parallel section of $D$ along the curve $\lambda_1$.
	Since $\dot\lambda_1(s)=-\xi_{\lambda_1(s)}$, we obtain 
	$$\nabla_{\dot\lambda_1(s)}^W\left(e^{sC_{\lambda_1(s)}}X(s) \right)=\left(\nabla_{\dot\lambda_1(s)}^{\boldsymbol{\tau}}-C_{\lambda_1(s)}\right)\left(e^{sC_{\lambda_1(s)}}X(s) \right)=0.$$
	This means that $$\tau^W_{\lambda_1} X(0)=e^{rC_x}X(r)=
	e^{rC_x}\tau^{\boldsymbol{\tau}}_{\lambda_1} X(0),$$
	and $$\tau^W_{\lambda_1} =e^{rC_x}\circ\tau^{\boldsymbol{\tau}}_{\lambda_1}.$$
	Moving to inverse maps, we obtain
	$$\tau^W_\lambda= 
	\tau^{\boldsymbol{\tau}}_\lambda\circ e^{-rC_x}.$$
	\qed

	\begin{lem}\label{lemxi2}\cite[Lemma 3.17]{GLL25} Let $\mu:[a,b]\to M$ be an arbitrary curve. Then there exists a horizontal curve $\bar\mu:[a,b]\to M$ such that
		$$\bar\mu(a)=\mu(a),\quad \bar\mu(b)=\mu(b),\quad \tau^g_{\bar\mu}=\tau^W_\mu.$$ 
	\end{lem}
	
	Let $\varphi_s$ denote the local flow of the vector field $\xi$. 
	
	\begin{lem}\label{lemxi3} Let $\mu:[a,b]\to M$ be a horizontal curve
		with $\mu(a)=x$, $\mu(b)=y$. Then for any $s\in\mathbb{R}$ there exists a horizontal curve $\bar\mu:[a,b]\to M$ 
		with $\bar\mu(a)=x$, $\bar\mu(b)=y$, such that
		\begin{equation}\label{taumusCx}\tau^g_\mu\circ e^{sC_x}=e^{sC_y}\circ \tau^g_{\bar\mu}.\end{equation}
	\end{lem}
	
	{\bf Proof.} Since $\mu([a,b])$ is compact, there exists such a positive $r\in\mathbb{R}$, that the map  
	$$\psi:[0,r]\times [a,b]\to M, \quad \psi(s,t)=\varphi_{s}\mu(t)$$
	is well-defined. A similar map was constructed in the proof of proposition 3.9 from \cite{GLL25}. Consider the curves 
	$$\check\mu(t)=\psi(r,t),\quad\lambda_1(s)=\psi(s,a),\quad \lambda_2(s)=\psi(s,b).$$
	It is obvious that the curve $\check\mu(t)$ is horizontal, and the curves $\lambda_1$ and $\lambda_2$ are integral curves of the vector field $\xi$.
	Consider the loop $$\mu*\lambda_2*\check\mu^{-1}*\lambda_1^{-1}.$$ The reasoning from the proof of proposition 3.9 in \cite{GLL25} shows that the parallel transport along this loop relative to the connection $\nabla^{\boldsymbol{\tau}}$ is an identity transformation. It follows from here that  
	$$\tau^{\boldsymbol{\tau}}_{\lambda_2}\circ\tau^{\boldsymbol{\tau}}_\mu=\tau^{\boldsymbol{\tau}}_{\check\mu}\circ\tau^{\boldsymbol{\tau}}_{\lambda_1}.$$
	By lemma \ref{lemxi1} we have
	$$\tau^W_{\lambda_2}\circ e^{rC_y}\circ\tau^g_\mu=\tau^g_{\check\mu}\circ \tau^W_{\lambda_1}\circ e^{rC_x}.$$
	According to lemma \ref{lemxi2}, there exists a horizontal curve $\bar\mu$, such that 
	$$\bar\mu(a)=\mu(a),\quad \bar\mu(b)=\mu(b),\quad \tau^g_{\bar\mu}=(\tau^W_{\lambda_2})^{-1}\circ\tau^g_{\check\mu}\circ \tau^W_{\lambda_1}.$$ 
	Consequently, $$e^{rC_y}\circ\tau^g_\mu=\tau^g_{\bar\mu}\circ  e^{rC_x}.$$
	This shows that the statement of the lemma is valid for all $s\in[0,r]$. 
	
	Let $H^{\boldsymbol{\tau}}_{x,y}$ be the set consisting of parallel transports relative to the connection $\nabla^{\boldsymbol{\tau}}$ along curves from $x$ to $y$. We fix a horizontal curve $\gamma$ from $x$ to $y$. Consider the map $$\operatorname{Hol}_x(\nabla^{\boldsymbol{\tau}})\to H^{\boldsymbol{\tau}}_{x,y},\quad \tau^{\boldsymbol{\tau}}_\mu\mapsto \tau^{\boldsymbol{\tau}}_\gamma\circ \tau^{\boldsymbol{\tau}}_\mu.$$ It is obvious that this map is a bijection. It defines a smooth manifold structure on the set $H^{\boldsymbol{\tau}}_{x,y}$, such that $H^{\boldsymbol{\tau}}_{x,y}$ is diffeomorphic to the manifold $\operatorname{Hol}_x(\nabla^{\boldsymbol{\tau}})$. 
	Similarly, let 
	$H^g_{x,y}\subset H^{\boldsymbol{\tau}}_{x,y}$   be a subset consisting of parallel transports relative to the connection $\nabla^g$ along horizontal curves from $x$ to $y$. We see that $H^g_{x,y}\subset H^{\boldsymbol{\tau}}_{x,y}$ is a submanifold, diffeomorphic to the manifold $\operatorname{Hol}_x(\nabla^g)$.
	
	The map
	$$(s,\tau^{\boldsymbol{\tau}}_\gamma) \mapsto  
	e^{sC_y}\circ\tau^{\boldsymbol{\tau}}_\gamma\circ e^{-sC_x}$$ defines an action of the Lie group $\mathbb{R}$ on the manifold $H^{\boldsymbol{\tau}}_{x,y}$. Let $\zeta$   be the induced vector field on $H^{\boldsymbol{\tau}}_{x,y}$. Since the statement of the lemma is valid for $s\in[0,r]$, the vector field $\zeta$ is tangent to $H^g_{x,y}$ at each point of $H^g_{x,y}$. This means that the action of the group $\mathbb{R}$ preserves the submanifold $H^g_{x,y}$, i. e. the statement of the lemma is valid for all $s\in\mathbb{R}$. \qed
	
	The rest of the proof is similar to the proof of theorem 3.14 from \cite{GLL25}.
	As in \cite{GalHolK,GLL25}, let us define the function
	$$f(t)=-\int_a^t\theta(\dot\mu(s))d s,\quad  t\in[a,b].$$
	Let $\mu:[a,b]\to M$   be a curve. It is clear that there exists $r\in(a,b]$, such that the curve
	$$\tilde\mu(t)=\varphi_{f(t)}\mu(t),\quad t\in [a,r]$$ is well-defined. According to \cite[Lemma 1]{GalHolK}, the curve $\mu(t)$ is horizontal.
	
	Let $\gamma:[a,b]\to M$   be an arbitrary loop at point $x$.
	Since $\gamma([a,b])\subset M$ is compact, there exist
	numbers $$a_0=a< a_1<\dots< a_{k-1}<a_k=b$$ such that
	for each restriction
	$$\gamma_i=\gamma|_{[a_{i-1},a_i]}$$ the horizontal curve $\tilde \gamma_i$ is well-defined. By definition of $\tilde\gamma_i$, there exists an integral curve $\lambda_i(s)$, $s\in [0,r_i]$ of the vector field $\xi$, connecting the final points of the curve $\tilde\gamma_i$ with the final point of the curve $\gamma_i$, and in this case $$r_i= \int_{\gamma_i}\theta.$$
	The reasoning from the proof of proposition 3.9 in \cite{GLL25} entails
	$$\tau^{\boldsymbol{\tau}}_{\gamma_i}=\tau^{\boldsymbol{\tau}}_{\lambda_i}\circ \tau_{\tilde\gamma_i}.$$
	It follows from lemmas \ref{lemxi1} and \ref{lemxi2} that there exists a horizontal curve $\bar\mu_i$, connecting points $\gamma(a_{i-1})$ and $\gamma(a_i)$, such that 
	$$\tau^{\boldsymbol{\tau}}_{\gamma_i}=e^{r_iC_{\gamma(a_i)}}\circ\tau_{\bar\mu_i}.$$
	The equality
	$$\tau^{\boldsymbol{\tau}}_\gamma=\tau^{\boldsymbol{\tau}}_{\lambda_k}\circ\cdots\circ\tau^{\boldsymbol{\tau}}_{\lambda_1},$$
	lemma \ref{lemxi3} and reasoning from the proof of theorem 3.14 from \cite{GLL25} entail the existence of a horizontal loop $\mu$ at point $x$, such that 
	$$\tau^{\boldsymbol{\tau}}_\gamma=\exp\left(\left(\int_\gamma\theta\right)\cdot {C_x}\right)\circ\tau_{\mu}.$$
	This, together with lemma \ref{lemxi3}, easily shows that the subgroup $\operatorname{Hol}_x(\nabla)\subset\operatorname{Hol}_x(\nabla^{\boldsymbol{\tau}})$ is normal. Finally, as in \cite{GLL25}, the map
	$$\lambda:\mathbb{R}_+\to\operatorname{Hol}_x(\nabla^{0})/\operatorname{Hol}_x(\nabla),$$
	$$\lambda:\exp\left(\int_\gamma\theta\right)\mapsto
	\exp\left(\left(\int_\gamma\theta\right)\cdot {C_x}\right)\cdot \operatorname{Hol}_x(\nabla)$$
	is well-defined and surjective. \qed
	
	\section{The case of K-contact sub-Lorentzian spaces}\label{SecLor}
	
	First of all, let us recall the classification of holonomy algebras of Lorentzian manifolds.
	We will use notations and exposition from \cite{Gal15}. The holonomy algebra of a Lorentzian manifold of dimension $k+2$ is contained in the Lie algebra $\mathfrak{so}(1,k+1)$. 
	The only irreducible Lorentzian holonomy algebra is $\mathfrak{so}(1,k+1)$.
	
	A subalgebra $\mathfrak{g}\subset \mathfrak{so}(1,k+1)$ is called weakly irreducible if it does not preserve any proper non-degenerate vector subspace in $\mathbb{R}^{1,k+1}$. If $\mathfrak{g}$ is weakly irreducible, but not irreducible, then $\mathfrak{g}$ preserves an isotropic line in $\mathbb{R}^{1,k+1}$.
	Let us fix a Witt basis $p,e_1,\dots,e_k,q$ of the space~$\mathbb{R}^{1,k+1}$.
	
	Let $\mathfrak{so}(1,k+1)_{\mathbb{R} p}$ denote the subalgebra in~$\mathfrak{so}(1,k+1)$, preserving the isotropic line~$\mathbb{R} p$. The Lie algebra $\mathfrak{so}(1,k+1)_{\mathbb{R} p}$ can be identified with the following matrix Lie algebra:
	$$
	\mathfrak{so}(1,k+1)_{\mathbb{R} p}=\left\{\begin{pmatrix} a &X^t
	& 0
	\\
	0 & A & -X
	\\
	0 & 0 & -a
	\end{pmatrix}\biggm| a\in \mathbb{R},\ X\in \mathbb{R}^k,\
	A \in \mathfrak{so}(k)\right\}.
	$$
	We identify the above-mentioned matrix with a triple $(a,A,X)$. This defines subalgebras~$\mathbb{R}$, $\mathfrak{so}(k)$, $\mathbb{R}^k$ in~$\mathfrak{so}(1,k+1)_{\mathbb{R}p}$. It is true that~$\mathbb{R}$ commutes with~$\mathfrak{so}(k)$, and~$\mathbb{R}^k$ is an ideal; it holds
	$$
	[(a,A,0),(0,0,X)]=(0,0,aX+AX).
	$$
	This gives the decomposition
	$$ \mathfrak{so}(1,k+1)_{\mathbb{R} p}=
	(\mathbb{R}\oplus\mathfrak{so}(k))\ltimes\mathbb{R}^k.
	$$
	Weakly irreducible, but not irreducible holonomy algebras of Lorentzian manifolds are exhausted by the following subalgebras in $ \mathfrak{so}(1,k+1)_{\mathbb{R} p}$ \cite{BBI,Gal06,Leistner}:
	
	{\rm\textbf{Type~1}}:
	$$
	\mathfrak{g}^{1,\mathfrak{h}}=(\mathbb{R}\oplus\mathfrak{h})\ltimes
	\mathbb{R}^k,
	$$
	
	{\rm\textbf{Type~2}}:
	$$
	\mathfrak{g}^{2,\mathfrak{h}}=\mathfrak{h}\ltimes\mathbb{R}^k;
	$$
	
	{\rm\textbf{Type~3}}:
	$$\mathfrak{g}^{3,\mathfrak{h},\varphi}=\{(\varphi(A),A,0)\mid
	A\in\mathfrak{h}\}\ltimes\mathbb{R}^k$$
	
	{\rm\textbf{Type~4}}:
	$$
	\mathfrak{g}^{4,\mathfrak{h},l,\psi}=\{(0,A,X+\psi(A))\mid
	A\in\mathfrak{h},\ X\in \mathbb{R}^m\}.$$
	
	Here $\mathfrak{h}\subset\mathfrak{so}(k)$   is a Riemannian holonomy algebra.
	In case of type 3, $\mathfrak{z}(\mathfrak{h})\ne\{0\}$, and
	$\varphi\colon\mathfrak{h}\to\mathbb{R}$   is a non-zero linear map
	with the property $\varphi\big|_{\mathfrak{h}'}=0$.
	In case of type 4 there exists an orthogonal decomposition
	$$\mathbb{R}^k=\mathbb{R}^l\oplus\mathbb{R}^{k-l}$$ such that
	$\mathfrak{h}\subset\mathfrak{so}(l)$,
	$\dim\mathfrak{z}(\mathfrak{h})\geqslant k-l$, and
	$\psi\colon\mathfrak{h}\to \mathbb{R}^{k-l}$   is a surjective linear map with the property $\psi\big|_{\mathfrak{h}'}=0$.
	
	The subalgebra $\mathfrak{h}\subset\mathfrak{so}(k)$, matched above to the weakly irreducible holonomy algebra $\mathfrak{g}\subset \mathfrak{so}(1,k+1)_{\mathbb{R} p}$, is called the {\bf orthogonal part} of the Lie algebra~$\mathfrak{g}$.
	For a Riemannian holonomy algebra $\mathfrak{h}\subset\mathfrak{so}(k)$ there exists an orthogonal decomposition
	\begin{equation} \label{eq4.4}
	\mathbb{R}^{k}=\mathbb{R}^{k_0}\oplus\mathbb{R}^{k_1}\oplus\cdots\oplus\mathbb{R}^{k_r}
	\end{equation}
	and the corresponding decomposition of $\mathfrak{h}$ into a direct sum of ideals
	\begin{equation}
	\label{eq4.5}
	\mathfrak{h}=\{0\}\oplus\mathfrak{h}_1\oplus\cdots\oplus\mathfrak{h}_r
	\end{equation}
	such that $\mathfrak{h}_i(\mathbb{R}^{k_j})=0$ at $i\ne j$, $\mathfrak{h}_i\subset\mathfrak{so}(k_i)$, and the representation of $\mathfrak{h}_i$ in $\mathbb{R}^{k_i}$ is irreducible.
	
	\begin{theorem}\label{ThLor} Let $\mathfrak{g}\subset\mathfrak{so}(1,k+1)_{\mathbb{R} p}$ be a weakly irreducible Lorentzian holonomy algebra. Then all ideals $I$ of codimension one in $\mathfrak{g}$ are contained in the following list:
		
		\begin{itemize}	
			\item[1.] $\mathfrak{g}=\mathfrak{g}^{1,\mathfrak{h}}$ and
			\begin{itemize}
				\item[1.1.] $I=\mathfrak{h}\ltimes \mathbb{R}^k$;
				\item[1.2.] $I=(\mathbb{R}\oplus I_1)\ltimes \mathbb{R}^k$, where $I_1\subset\mathfrak{h}$   is an ideal of codimension one;
				\item[1.3.] $I=\mathfrak{g}^{3,\mathfrak{h},\varphi}$ for some non-zero linear map $\varphi:\mathfrak{h}\to\mathbb{R}$;
			\end{itemize}
			\item[2.] $\mathfrak{g}=\mathfrak{g}^{2,\mathfrak{h}}$ and
			\begin{itemize}
				\item[2.1.] $I=\mathfrak{h}\ltimes \mathbb{R}^{k-1}$, where $\mathfrak{h}\subset\mathfrak{so}(k-1)$;
				\item[2.2.] $I=I_1\ltimes \mathbb{R}^k$, where $I_1\subset\mathfrak{h}$   is an ideal of codimension one;
				\item[2.3.] $I=\mathfrak{g}^{4,\mathfrak{h},\varphi}=\{A+\psi(A)|A\in\mathfrak{h}\}\ltimes\mathbb{R}^{k-1}$, where $\mathfrak{h}\subset\mathfrak{so}(k-1)$, and $\psi:\mathfrak{h}\to\mathbb{R}$   is a non-zero linear map;
			\end{itemize}
			\item[3.] $\mathfrak{g}=\mathfrak{g}^{3,\mathfrak{h},\varphi}$ and
			\begin{itemize}
				\item[3.1.] $I=\ker\varphi\ltimes \mathbb{R}^k$; 
				\item[3.2.] $I=\mathfrak{g}^{3,I_1,\varphi_1}$, where $I_1\subset\mathfrak{h}$   is an ideal of codimension one such that $\varphi_1=\varphi|_{I_1}$ is not zero;
			\end{itemize}
			\item[4.] $\mathfrak{g}=\mathfrak{g}^{4,\mathfrak{h},l,\psi}$
			\begin{itemize}
				\item[4.1.] $I=\{A+\psi(A)|A\in\mathfrak{h}\}\ltimes \mathbb{R}^{l-1}$, where $\mathfrak{h}\subset\mathfrak{so}(l-1)$;
				\item[4.2.] $I=\{A+\psi(A)|A\in I_1\}\ltimes \mathbb{R}^{l}$, where $I_1\subset\mathfrak{h}$   is an ideal of codimension one;
				\item[4.3.] $I=\{A+\psi_1(A)|A\in\mathfrak{h}\} \ltimes \mathbb{R}^{l-1}$, where $\psi_1:\mathfrak{h}\to\mathbb{R}^{k-l+1}$   is a surjective linear map such that $\psi=\operatorname{pr}_{\mathbb{R}^{k-l}}\circ\psi_1$.
			\end{itemize}
		\end{itemize}
	\end{theorem}
	
	{\bf Proof.} We will consider 4 cases according to the type of $\mathfrak{g}$.
	
	{\bf Case 1.} Suppose that $\mathfrak{g}=\mathfrak{g}^{1,\mathfrak{h}}$. First, let's assume that the projection of $I$ on $\mathbb{R}\subset\mathfrak{g}$ is trivial, i.e. $I\subset\mathfrak{h}\ltimes\mathbb{R}^k$. Since the codimension of $I$ in $\mathfrak{g}$ is one, we obtain $I=\mathfrak{h}\ltimes\mathbb{R}^k$, which corresponds to case 1.1.
	
	Suppose that the projection of $I$ on $\mathbb{R}\subset\mathfrak{g}$ is non-trivial. Then there exists an element $(a,A,X)\in I$ with $a\neq 0$. Let $Y\in\mathbb{R}^k\subset \mathfrak{g}$. Since $I\subset\mathfrak{g}$   is an ideal, we obtain that
	$$[(a,A,X),Y]=(0,0,aY+AY)\in I.$$
	Since $A\in\mathfrak{h}\subset\mathfrak{so}(k)$ has no non-zero real eigenvalues, the map
	$$Y\in\mathbb{R}^k\mapsto  aY+AY\in\mathbb{R}^k$$ is an isomorphism. This shows that $\mathbb{R}^k\subset I$. We see that $I=I_0\ltimes\mathbb{R}^k$, where $I_0\subset\mathbb{R}\oplus\mathfrak{h}$   is an ideal of codimension one. It is clear that $I$ corresponds to case 1.2 or 1.3.
	
	{\bf Case 3.} Suppose that $\mathfrak{g}=\mathfrak{g}^{3,\mathfrak{h},\varphi}$. The reasoning is similar to the previous case: if the projection of $I$ on $\mathbb{R}\subset\mathfrak{so}(1,k+1)$ is trivial, then $I$ corresponds to case 3.1; if the projection of $I$ on $\mathbb{R}\subset\mathfrak{so}(1,k+1)$ is non-trivial, then $I$ corresponds to case 3.2.
	
	{\bf Case 2.} For $\mathfrak{g}$ of type 2 and 4 it is convenient to consider on $\mathbb{R}\oplus\mathfrak{so}(k)\oplus\mathbb{R}^k$ a scalar product, in which this direct sum is orthogonal, and represent $I$ as $\alpha^\bot$ for some non-zero $\alpha\in \mathfrak{g}$. 
	
	Consider the above-mentioned decomposition \eqref{eq4.4}. Since $k_1,\dots,k_r\geqslant 2$ and $I\subset\mathfrak{g}$ has codimension one, it holds $\operatorname{pr}_{\mathbb{R}^{k_i}}I\neq\{0\}$ for $i=1,\dots,r$. Since for each $i$ the algebra $\mathfrak{h}_i\subset\mathfrak{so}(k_i)$ is irreducible, and $I\subset\mathfrak{g}$ is an ideal, we conclude that
	$$\mathbb{R}^{k_1}\oplus\cdots\oplus\mathbb{R}^{k_r}\subset I.$$ We see that $\alpha\in\mathfrak{h}\oplus\mathbb{R}^{k_0}$. If $\alpha\in \mathbb{R}^{k_0}$, then $I$ corresponds to case 2.1. If $\alpha\in \mathfrak{h}$, then $I$ corresponds to case 2.2. If the projections of $\alpha$ on $\mathfrak{h}$ and on $\mathbb{R}^{k_0}$ are both non-zero, then $I$ corresponds to case 2.3.
	
	{\bf Case 4.} Suppose that $\mathfrak{g}=\mathfrak{g}^{4,\mathfrak{h},\psi}$. Since $\mathfrak{h}\subset\mathfrak{so}(l)$, the decomposition \eqref{eq4.4} has the form 
	$$\mathbb{R}^k=\mathbb{R}^{k-l}\oplus\mathbb{R}^{l_0}\oplus\mathbb{R}^{l_1}\oplus\cdots\oplus\mathbb{R}^{l_r},$$
	$$\mathfrak{h}=\mathfrak{h}_1\oplus\cdots\oplus\mathfrak{h}_r,$$
	where $\mathfrak{h}_i\subset\mathfrak{so}(l_i)$   are irreducible subalgebras. As in case 2, it can be proven that
	$$\mathbb{R}^{l_1}\oplus\cdots\oplus\mathbb{R}^{l_r}\subset I.$$
	We see that $I=\alpha^\bot$, where $$\alpha\in\{(0,A,\psi(A))|A\in\mathfrak{h}\}\ltimes\mathbb{R}^{l_0}.$$
	If $\alpha\in\mathbb{R}^{l_0}$, then $I$ corresponds to case 4.1. If $\alpha\in\{(0,A,\psi(A))|A\in\mathfrak{h}\}$, then $I$ corresponds to case 4.2. In the last obvious possibility for $\alpha$ the ideal $I$ corresponds to case 4.3. \qed
	
	{\bf Example 1.} Consider the space $\mathbb{R}^{1,2s+1}$ with coordinates $v,x^1,x^2,\dots,x^{2s},u$ and metric
	$$h=2dvdu+\sum_{i=1}^{2s} (dx^i)^2+H(du)^2,\quad H=\sum_{i=1}^{2s} (x^i)^2.$$
	This is an example of a symmetric Cahen-Wallach space. See, for example, \cite{Gal15} for formulas of Christoffel symbols, curvature and holonomy of such metrics.
	The holonomy algebra $\mathfrak{hol}(\nabla^h)$ of this metric is isomorphic to $\mathbb{R}^{2s}\subset\mathfrak{so}(1,2s+1)$. The vector field $\partial_v$ is parallel:
	$$\nabla^h\partial_v=0.$$
	It holds
	$$\nabla^h_X\partial_{x^i}\in\left<\partial_v\right>.$$
	The curvature tensor of this space satisfies the conditions $$R^h(\partial_{x^i},\partial_u)=\partial_{x^i}\wedge \partial_v,\quad R^h(\partial_{x^i},\partial_{x^j})=R^h(\partial_v,\cdot)=0.$$
	Note that the bivectors $\partial_{x^i}\wedge \partial_v$ are parallel.
	The holonomy algebra at an arbitrary point is spanned by the bivectors $\partial_{x^i}\wedge \partial_v$ at this point. 
	Let $m=s+1$ and $n=2m+1$.
	Let $$M=\mathbb{R}\times \mathbb{R}^{1,2s+1}$$
	and $$\theta=dt+\sum_{i=1}^{s-1} x^{2i-1}dx^{2i}  + v dx^{2s-1}+u dx^{2s}.$$
	It is obvious that $\theta$   is a contact form on $M$. As in \cite{Kokin1}, the fibers of the contact distribution $D=\ker\theta$ can be identified with tangent spaces to $\mathbb{R}^{1,2s+1}$, and then $h$ defines a K-contact sub-Lorentzian metric $g$ on $M$. Moreover, it holds $$\mathfrak{hol}(\nabla^{\boldsymbol{\tau}})\cong\mathfrak{hol}(\nabla^h).$$
	The distribution $D$ is spanned by the vector fields
	$$V=\partial_v,\quad X_{2i-1}=\partial_{x^{2i-1}},\quad X_{2i}=\partial_{x^{2i}}-x^{2i-1}\partial_t,\quad i=1,\dots,s-1,$$
	$$ X_{2s-1}=\partial_{x^{2s-1}}-v\partial_t,\quad X_{2s}=\partial_{x^{2s}}-u\partial_t,\quad U=\partial_u.$$
	Dual 1-forms have the form:
	$$\alpha=dv,\quad \gamma^{2i-1}=dx^{2i-1},\quad \gamma^{2i}=d{x^{2i}},\quad i=1,\dots,s-1,$$
	$$\gamma^{2s-1}=d{x^{2s-1}},\quad \gamma^{2s}=d{x^{2s}},\quad \beta=du.$$
	It holds $$d\theta=\sum_{i=1}^{s-1} \gamma^{2i-1}\wedge \gamma^{2i}+\alpha\wedge \gamma^{2s-1}+\beta\wedge X^{2s}.$$
	For the bivector $(d\theta)^{-1}$, defined above, we have
	$$(d\theta)^{-1}=-2\left(\sum_{i=1}^{s-1} X_{2i-1}\wedge X_{2i}+V\wedge X_{2s-1}+U\wedge X_{2s}\right).$$
	
	According to the Ambrose-Singer theorem for the horizontal holonomy algebra $\mathfrak{hol}_x(\nabla^g)$ at point $x\in M$ \cite{CGJK,FGR}, $\mathfrak{hol}_x(\nabla^g)$ is spanned by the elements
	$$(\tau_\gamma^g)^{-1}\circ R^g_y(B)\circ\tau_\gamma^g,$$
	where $\gamma$   is an arbitrary horizontal curve starting at $x$, $y$   final point of $\gamma$, and $B$   a bivector from $D_y$ such that $$d\theta(B)=0.$$
	Non-zero values of $R^g$ are
	$$R^g(X_i,U)=X_i\wedge V.$$
	It is clear that $d\theta$ annuls the bivectors
	$$U\wedge X_1,\dots, U\wedge X_{2s-1}.$$
	On the other hand, $$d\theta(U\wedge X_{2s})=2.$$
	The bivectors $X_i\wedge V$, $i=1,\dots,2s$, are $\nabla^g$-parallel.
	We conclude that
	$$\mathfrak{hol}_x(\nabla^g)=\left<X_1\wedge V,\dots,X_{2s-1}\wedge V\right>_x.$$
	Thus, $$ \mathfrak{hol}_x(\nabla^g)\cong\mathbb{R}^{2s-1}\subset\mathfrak{hol}_x(\nabla^{\boldsymbol{\tau}})\cong\mathbb{R}^{2s}$$
	has codimension one. This is a particular case of 2.1 from theorem \ref{ThLor}.
	
	{\bf Example 2.} Let $(N_0,h_0)$   be a locally symmetric Kähler manifold with an irreducible holonomy algebra $\mathfrak{h}\subset\mathfrak{u}(s)$. It holds $$\mathfrak{h}=\mathfrak{h}'\oplus\mathbb{R} J,\quad \mathfrak{h}'\subset\mathfrak{su}(s).$$ Suppose that there exists a 1-form $\theta_0$, such that $d\theta_0$ is a Kähler form on $(N_0,h_0)$. Consider the manifold
	$$N=\mathbb{R}\times N_0\times \mathbb{R}$$ with Lorentzian metric
	$$h=2dvdu+h_0+H(du)^2,$$
	where $H$   is a function on $N_0$. The holonomy algebra of $N$ is 
	$$\mathfrak{h}\ltimes\mathbb{R}^{2s}\subset\mathfrak{so}(1,2s+1)$$ for a sufficiently general function $H$.
	Let $$M=\mathbb{R}\times N$$   manifold with contact form $$\theta=dt+vdu+\theta_0.$$
	Again, as in \cite{Kokin1}, the fibers of the contact distribution $D=\ker\theta$ are identified with tangent spaces to $N$, and $h$ defines a K-contact sub-Lorentzian metric $g$ on $M$. It holds $$\mathfrak{hol}(\nabla^{\boldsymbol{\tau}})\cong\mathfrak{hol}(\nabla^h).$$
	The distribution $D$ is spanned by the vector fields
	$$V=\partial_v,\quad X_{i}=\partial_{x^{i}}-\theta_0(\partial_{x^{i}})\partial_t,\quad i=1,\dots,2s,\quad U=\partial_u-v\partial_t.$$
	Dual 1-forms have the form:
	$$\alpha=dv,\quad \gamma^{i}=dx^{i},\quad i=1,\dots,2s,
	\quad \beta=du.$$
	It holds $$d\theta=d\theta_0+\alpha\wedge\beta.$$
	Let $J$   be a complex structure on $(N_0,h_0)$. We consider $J$ as an endomorphism field of $D$.
	As it is known, $d\theta_0(J)=2s$. This entails
	$$d\theta(J)=d\theta_0(J)=2s.$$ In \cite{GLL25} it is shown that
	$$d\theta((d\theta)^{-1})=-4m.$$
	We conclude that 
	$$(d\theta)^{-1}=-2J-2V\wedge U.$$
	It is clear that $$\mathfrak{hol}_x(\nabla^g)\subset\mathfrak{hol}_x(\nabla^{\boldsymbol{\tau}})=\mathfrak{h}\ltimes\mathbb{R}^{2s}.$$
	As shown in \cite{Leistner1}, the orthogonal part $\mathfrak{h}$ of the holonomy algebra of the manifold $(N,h)$ coincides with the holonomy algebra of the induced connection on the so-called screen bundle
	$$\mathcal{E}=\left<\partial_v\right>^\bot/\left<\partial_v\right>.$$ The curvature tensor $R^\mathcal{E}$ of this connection satisfies the conditions 
	$$R^\mathcal{E}(\partial_{x^i},\partial_{x^j})=R_0(\partial_{x^i},\partial_{x^j}),\quad  R^\mathcal{E}(\partial_v,\cdot)=0,\quad R^\mathcal{E}(\partial_u,\cdot)=0,$$ where $R_0$   is the curvature tensor of the metric $h_0$.
	Similarly, the orthogonal part of $\mathfrak{hol}_x(\nabla^g)$ coincides with the holonomy algebra $\mathfrak{hol}_x(\nabla^{\mathcal{E}_1})$ of the induced horizontal connection on the bundle 
	$$\mathcal{E}_1=\left<V\right>^\bot/\left<V\right>.$$
	As in the previous example, by the Ambrose-Singer theorem the horizontal holonomy $\mathfrak{hol}_x(\nabla^{\mathcal{E}_1})$ is spanned by the elements
	$$(\tau_\gamma^{\mathcal{E}_1})^{-1}\circ R^{\mathcal{E}_1}_y(B)\circ\tau_\gamma^{\mathcal{E}_1},$$
	where $\gamma$   is an arbitrary horizontal curve starting at $x$, $y$   final point of $\gamma$, and $B$   a bivector from $D_y$ such that $$d\theta(B)=0.$$
	Since $R^{\mathcal{E}_1}$ is completely defined through $R_0$, we have $$R^{\mathcal{E}_1}_y(B)\in\mathfrak{su}({\mathcal{E}_1}_y)\quad \text{if } d\theta(B)=0.$$ Thus, $\mathfrak{hol}_x(\nabla^{\mathcal{E}_1})=\mathfrak{h}'$ and 
	$$\mathfrak{hol}_x(\nabla^g)=\mathfrak{h}'\ltimes \mathbb{R}^{2s}.$$
	We are in case 2.2 from theorem \ref{ThLor}.

\end{document}